\documentclass[12pt]{article}
\usepackage{amssymb}
\usepackage{amsfonts}
\usepackage{amsmath}
\usepackage{amsbsy}
\newcommand{\be}{\begin{equation}}
\newcommand{\ee}{\end{equation}}
\newcommand{\bea}{\begin{eqnarray}}
\newcommand{\eea}{\end{eqnarray}}
\newcommand{\ba}{\begin{array}}
\newcommand{\ea}{\end{array}}

\newcommand{\bc}{\begin{center}}
\newcommand{\ec}{\end{center}}
\newcommand{\ben}{\begin{enumerate}}
\newcommand{\een}{\end{enumerate}}
\newcommand{\bfi}{\begin{figure}}
\newcommand{\efi}{\end{figure}}

\newcommand{\bq}{\begin{quote}}
\newcommand{\eq}{\end{quote}}
\newcommand{\bqu}{\begin{quotation}}
\newcommand{\equ}{\end{quotation}}
\newenvironment{emphit}{\begin{itemize}}{\end{itemize}}
\newcommand{\bemp}{\begin{emphit}}
\newcommand{\eemp}{\end{emphit}}

\newcommand{\bt}{\begin{tabular}}
\newcommand{\et}{\end{tabular}}

\newtheorem{myth}{Theorem}[section]

\newtheorem{mydef}{Definition}[section]
\newtheorem{myrem}{Remark}[section]
\newtheorem{myexam}{Example}[section]

\begin{document}
\date{}
\title{A Pair Of Smarandachely Isotopic Quasigroups And Loops Of The Same Variety
}
\author{T\`em\'it\d{\'o}p\d{\'e} Gb\d{\'o}l\'ah\`an Ja\'iy\'e\d ol\'a\thanks{All correspondence to be addressed to this author.} \\
Department of Mathematics,\\
Obafemi Awolowo University, Ile-Ife 220005, Nigeria.\\
jaiyeolatemitope@yahoo.com,~tjayeola@oauife.edu.ng} \maketitle
\begin{abstract}
The isotopic invariance or universality of types and varieties of
quasigroups and loops described by one or more equivalent identities
has been of interest to researchers in loop theory in the recent
past. A variety of quasigroups(loops) that are not universal have
been found to be isotopic invariant relative to a special type of
isotopism or the other. Presently, there are two outstanding open
problems on universality of loops: semi automorphic inverse property
loops(1999) and Osborn loops(2005). Smarandache
isotopism(S-isotopism) was originally introduced by Vasantha
Kandasamy in 2002. But in this work, the concept is re-restructured
in order to make it more explorable. As a result of this, the theory
of Smarandache isotopy inherits the open problems as highlighted
above for isotopy. In this short note, the question 'Under what type
of S-isotopism will a pair of S-quasigroups(S-loops) form any
variety?' is answered by presenting a pair of specially S-isotopic
S-quasigroups(loops) that both belong to the same variety of
S-quasigroups(S-loops). This is important because pairs of specially
S-isotopic S-quasigroups(e.g Smarandache cross inverse property
quasigroups) that are of the same variety are useful for
applications(e.g cryptography).
\end{abstract}
\newpage
\section{Introduction}
\subsection{Isotopy Theory Of Quasigroups And Loops}
The isotopic invariance of types and varieties of quasigroups and
loops described by one or more equivalent identities, especially
those that fall in the class of Bol-Moufang type loops as first
named by Fenyves \cite{phd56} and \cite{phd50} in the 1960s  and
later on in this $21^\textrm{st}$ century by Phillips and Vojt\v
echovsk\'y \cite{phd9}, \cite{phd61} and \cite{phd124} have been of
interest to researchers in loop theory in the recent past. Among
such is Etta Falconer's Ph.D \cite{phd159} and her paper
\cite{phd160} which investigated isotopy invariants in quasigroups.
Loops such as Bol loops, Moufang loops, central loops and extra
loops are the most popular loops of Bol-Moufang type whose isotopic
invariance have been considered. For more on loops and their
properties, readers should check \cite{phd3},
\cite{phd41},\cite{phd39}, \cite{phd49}, \cite{phd42} and
\cite{phd75}.

Bol-Moufang type of quasigroups(loops) are not the only
quasigroups(loops) that are isomorphic invariant and whose
universality have been considered. Some others are flexible loops,
F-quasigroups, totally symmetric quasigroups(TSQ), distributive
quasigroups, weak inverse property loops(WIPLs), cross inverse
property loops(CIPLs), semi-automorphic inverse property
loops(SAIPLs) and inverse property loops(IPLs). As shown in Bruck
\cite{phd3}, a left(right) inverse property loop is universal if and
only if it is a left(right) Bol loop, so an IPL is universal if and
only if it is a Moufang loop. Ja\'iy\'e\d ol\'a \cite{phd143}
investigated the universality of central loops. Recently, Michael
Kinyon et. al. \cite{phd95}, \cite{phd118}, \cite{phd119} solved the
Belousov problem concerning the universality of F-quasigroup which
has been open since 1967. The universality of WIPLs and CIPLs have
been addressed by OSborn \cite{phd89} and Artzy \cite{phd30}
respectively while the universality of elasticity(flexibility) was
studied by Syrbu \cite{phd97}. In 1970, Basarab \cite{phd146} later
continued the work of J. M. Osborn of 1961 on universal WIPLs by
studying isotopes of WIPLs that are also WIPLs after he had studied
a class of WIPLs(\cite{phd149}) in 1967. The universality of SAIPLs
is still an open problem to be solved as stated by Michael Kinyon
during the LOOPs'99 conference. After the consideration of universal
AIPLs by Karklinsh and Klin \cite{phd150}, Basarab \cite{phd147}
obtained a sufficient condition for which a universal AIPL is a
G-loop. Although Basarab \cite{phd148} and \cite{phd170} considered
universal Osborn loops but the universality of Osborn loops was
raised as an open problem by Michael Kinyon \cite{phd33}. Up to the
present moment, the problem is still open.

Interestingly, Adeniran \cite{phd79} and Robinson \cite{phd85},
Oyebo and Adeniran \cite{phd141}, Chiboka and Solarin \cite{phd80},
Bruck \cite{phd82}, Bruck and Paige \cite{phd40}, Robinson
\cite{phd7}, Huthnance \cite{phd44} and Adeniran \cite{phd79} have
respectively studied the holomorphs of Bol loops, central loops,
conjugacy closed loops, inverse property loops, A-loops, extra
loops, weak inverse property loops, Osborn loops and Bruck loops.
Huthnance \cite{phd44} showed that if $(L,\cdot )$ is a loop with
holomorph $(H,\circ)$, $(L,\cdot )$ is a WIPL if and only if
$(H,\circ)$ is a WIPL. The holomorphs of an AIPL and a CIPL are yet
to be studied.

\subsection{Isotopy Theory Of Smarandache Quasigroups And Loops}
The study of Smarandache loops was initiated by W.B. Vasantha
Kandasamy in 2002. In her book \cite{phd75}, she defined a
Smarandache loop(S-loop) as a loop with at least a subloop which
forms a subgroup under the binary operation of the loop. In her
book, she introduced over 75 Smarandache concepts on loops. In her
first paper \cite{phd83}, she introduced Smarandache : left(right)
alternative loops, Bol loops, Moufang loops, and Bruck loops. But in
Ja\'iy\'e\d ol\'a \cite{sma1}, Smarandache : inverse property loops
(IPL), weak inverse property loops (WIPL), G-loops, conjugacy closed
loops (CC-loop), central loops, extra loops, A-loops, K-loops, Bruck
loops, Kikkawa loops, Burn loops and homogeneous loops were
introduced and studied relative to the holomorphs of loops. It is
particularly established that a loop is a Smarandache loop if and
only if its holomorph is a Smarandache loop. This statement was also
shown to be true for some weak Smarandache loops(inverse property,
weak inverse property) but false for others(conjugacy closed, Bol,
central, extra, Burn, A-, homogeneous) except if their holomorphs
are nuclear or central. The study of Smarandache quasigroups was
carried out in Ja\'iy\'e\d ol\'a \cite{sma2} after their
introduction in Muktibodh \cite{muk1} and \cite{muk2}. In
Ja\'iy\'e\d ol\'a \cite{sma3}, the universality of some Smarandache
loops of Bol-Moufang types was studied and neccessary and sufficient
conditions for their universality were established.

\paragraph{}
In this short note, the question 'Under what type of S-isotopism
will a pair of S-quasigroups(S-loops) form any variety?' is answered
by presenting a pair of specially S-isotopic S-quasigroups(loops)
that both belong to the same variety of S-quasigroups(S-loops). This
is important because pairs of specially S-isotopic S-quasigroups(e.g
Smarandache cross inverse property quasigroups) that are of the same
variety are useful for applications(e.g cryptography).

\section{Definitions and Notations}
\begin{mydef}\label{1:0}
Let $L$ be a non-empty set. Define a binary operation ($\cdot $) on
$L$ : If $x\cdot y\in L~\forall ~x, y\in L$, $(L, \cdot )$ is called
a groupoid. If the system of equations ; $a\cdot x=b$ and $y\cdot
a=b$ have unique solutions for $x$ and $y$ respectively, then $(L,
\cdot )$ is called a quasigroup. Furthermore, if there exists a
unique element $e\in L$ called the identity element such that
$\forall
~x\in L$, $x\cdot e=e\cdot x=x$, $(L, \cdot )$ is called a loop.\\

If there exists at least a non-empty and non-trivial subset $M$ of a
groupoid(quasigroup or semigroup or loop) $L$ such that $(M,\cdot )$
is a non-trivial subsemigroup(subgroup or subgroup or subgroup) of
$(L, \cdot )$, then $L$ is called a Smarandache:
groupoid(S-groupoid)$\Big($quasigroup(S-quasigroup) or
semigroup(S-semigroup) or loop(S-loop)$\Big)$ with Smarandache:
subsemigroup(S-subsemigroup)$\Big($subgroup(S-subgroup) or
subgroup(S-subgroup) or subgroup(S-subgroup)$\Big)$ $M$.

A quasigroup(loop) is called a Smarandache "certain"
quasigroup(loop) if it has at least a non-trivial
subquasigroup(subloop) with the "certain" property and the latter is
referred to as the Smarandache "certain" subquasigroup(subloop). For
example, a loop is called a Smarandache Bol-loop if it has at least
a non-trivial subloop that is a Bol-loop and the latter is referred
to as the Smarandache Bol-subloop. By an "initial S-quasigroup" $L$
with an "initial S-subquasigroup" $L'$, we mean $L$ and $L'$ are
pure quasigroups, i.e. they do not obey a "certain" property(not of
any variety).

Let $(G,\cdot )$ be a quasigroup(loop). The bijection $L_x :
G\rightarrow G$ defined as $yL_x=x\cdot y~\forall ~x, y\in G$ is
called a left translation(multiplication) of $G$ while the bijection
$R_x : G\rightarrow G$ defined as $yR_x=y\cdot x~\forall ~x, y\in G$
is called a right translation(multiplication) of $G$.

The set $SYM(L, \cdot )=SYM(L)$ of all bijections in a groupoid
$(L,\cdot )$ forms a group called the permutation(symmetric) group
of the groupoid $(L,\cdot )$. If $L$ is a S-groupoid with a
S-subsemigroup $H$, then the set $SSYM(L, \cdot )=SSYM(L)$ of all
bijections $A$ in $L$ such that $A~:~H\to H$ forms a group called
the Smarandache permutation(symmetric) group of the S-groupoid. In
fact, $SSYM(L)\le SYM(L)$.
\end{mydef}

\begin{mydef}\label{1:1}
If $(L, \cdot )$ and $(G, \circ )$ are two distinct groupoids, then
the triple $(U, V, W) : (L, \cdot )\rightarrow (G, \circ )$ such
that $U, V, W : L\rightarrow G$ are bijections is called an
isotopism if and only if
\begin{displaymath}
xU\circ yV=(x\cdot y)W~\forall ~x, y\in L.
\end{displaymath}
So we call $L$ and $G$ groupoid isotopes.

If $U=V=W$, then $U$ is called an isomorphism, hence we write $(L,
\cdot )\cong (G, \circ )$.

Now, if $(L, \cdot )$ and $(G, \circ )$ are S-groupoids with
S-subsemigroups $L'$ and $G'$ respectively such that $A~:~L'\to G'$,
where $A\in\{U,V,W\}$, then the isotopism $(U, V, W) : (L, \cdot
)\rightarrow (G, \circ )$ is called a Smarandache
isotopism(S-isotopism).

Thus, if $U=V=W$, then $U$ is called a Smarandache isomorphism,
hence we write $(L, \cdot )\succsim (G, \circ )$.

If $(L,\cdot )=(G,\circ )$, then the triple $\alpha =(U, V, W)$ of
bijections on $(L,\cdot )$ is called an autotopism of the
groupoid(quasigroup, loop) $(L,\cdot )$. Such triples form a group
$AUT(L,\cdot )$ called the autotopism group of $(L,\cdot )$.
Furthermore, if $U=V=W$, then $U$ is called an automorphism of the
groupoid(quasigroup, loop) $(L,\cdot )$. Such bijections form a
group $AUM(L,\cdot )$ called the automorphism group of $(L,\cdot )$.

Similarly, if $(L, \cdot )$ is an S-groupoid with S-subsemigroup
$L'$ such that $A\in\{U,V,W\}$ is a Smarandache permutation, then
the autotopism $(U, V, W)$ is called a Smarandache autotopism
(S-autotopism) and they form a group $SAUT(L,\cdot )$ which will be
called the Smarandache autotopism group of $(L, \cdot )$. Observe
that $SAUT(L,\cdot )\le AUT(L,\cdot )$.
\end{mydef}

\paragraph{Discussions}
To be more precise about the notion of S-isotopism in
Definition~\ref{1:1}, the following explanations are given. For a
given S-groupoid, the S-subsemigroup is arbitrary. But in the
proofs, we make use of one arbitrary S-subsemigroup for an
S-groupoid at a time for our arguments. Now, if $(L, \cdot )$ and
$(G, \circ )$ are S-isotopic S-groupoids with arbitrary
S-subsemigroups $L'$ and $G'$ respectively under the triple
$(U,V,W)$. In case the S-subsemigroup $L'$ of the S-groupoid $L$ is
replaced with another S-groupoid $L''$ of $L$(i.e a situation where
by $L$ has at least two S-subsemigroups), then under the same
S-isotopism $(U,V,W)$, the S-groupoid isotope $G$ has a second
S-subsemigroups $G''$. Hence, when studying the S-isotopism
$(U,V,W)$, it will be for the system
\begin{displaymath}
\{(L, \cdot ), (L', \cdot )\} \rightarrow \{(G, \circ ),(G', \circ
)\}~\textrm{or}~ \{(L, \cdot ), (L'', \cdot )\} \rightarrow \{(G,
\circ ),(G'', \circ )\}
\end{displaymath}
and not
\begin{displaymath}
\{(L, \cdot ), (L', \cdot )\} \rightarrow \{(G, \circ ),(G'', \circ
)\}~\textrm{or}~ \{(L, \cdot ), (L'', \cdot )\} \rightarrow \{(G,
\circ ),(G', \circ )\}.
\end{displaymath}
This is because $|L'|=|G'|$ and $|L''|=|G''|$ since $(L')A=G'$ and
$(L'')A=G''$ for all $A\in\{U,V,W\}$ while it is not compulsory that
$|L'|=|G''|$ and $|L''|=|G'|$. It is very easy to see from the
definition that the component transformations $U,V,W$ of isotopy
after restricting them to the S-subsemigroup or S-subgroup $L'$ are
bijections. Let $x_1,x_2\in L'$, then $x_1A=x_2A$ implies that
$x_1=x_2$ because $x_1,x_2\in L'$ implies $x_1,x_2\in L$, hence
$x_1A=x_2A$ in $L$ implies $x_1=x_2$. The mappings $A~:L'\to G'$ and
$A~:L-L'\to G-G'$ are bijections because $A~:L\to G$ is a bijection.
Our explanations above are illustrated with the following examples.

\begin{myexam}\label{exam:1}
The systems $(L,\cdot )$ and $(L,\ast)$, $L=\{0,1,2,3,4\}$ with the
multiplication tables below are S-quasigroups with S-subgroups
$(L',\cdot )$ and $(L'',\ast )$ respectively, $L'=\{0,1\}$ and
$L''=\{1,2\}$. $(L,\cdot )$ is taken from Example~2.2 of
\cite{muk2}. The triple $(U,V,W)$ such that
\begin{displaymath}
U=\left(\begin{array}{ccccc} 0 & 1 & 2 & 3 & 4 \\
1 & 2 & 3 & 4 & 0
\end{array}\right),~V=
\left(\begin{array}{ccccc}
0 & 1 & 2 & 3 & 4 \\
1 & 2 & 4 & 0 & 3
\end{array}\right)~\textrm{and}~W=
\left(\begin{array}{ccccc}
0 & 1 & 2 & 3 & 4 \\
1 & 2 & 0 & 4 & 3
\end{array}\right)
\end{displaymath}
are permutations on $L$, is an S-isotopism of $(L,\cdot )$ onto
$(L,\ast)$. Notice that $A(L')=L''$ for all $A\in\{U,V,W\}$
and $U, V, W : L'\rightarrow L''$ are all bijcetions.\\
\begin{center}
\begin{tabular}{|c||c|c|c|c|c|}
\hline
$\cdot $ & 0 & 1 & 2 & 3 & 4 \\
\hline \hline
0 & 0 & 1 & 3 & 4 & 2 \\
\hline
1 & 1 & 0 & 2 & 3 & 4 \\
\hline
2 & 3 & 4 & 1 & 2 & 0 \\
\hline
3 & 4 & 2 & 0 & 1 & 3 \\
\hline
4 & 2 & 3 & 4 & 0 & 1 \\
\hline
\end{tabular}
\qquad\qquad\qquad\qquad\qquad\qquad\qquad\qquad\begin{tabular}{|c||c|c|c|c|c|}
\hline
$\ast$ & 0 & 1 & 2 & 3 & 4 \\
\hline \hline
0 & 1 & 0 & 4 & 2 & 3 \\
\hline
1 & 3 & 1 & 2 & 0 & 4 \\
\hline
2 & 4 & 2 & 1 & 3 & 0 \\
\hline
3 & 0 & 4 & 3 & 1 & 2 \\
\hline
4 & 2 & 3 & 0 & 4 & 1 \\
\hline
\end{tabular}
\end{center}
\end{myexam}

\begin{myexam}\label{exam:2}
According Example~4.2.2 of \cite{van2}, the system
$(\mathbb{Z}_6,\times_6)$ i.e the set $L=\mathbb{Z}_6$ under
multiplication modulo $6$ is an S-semigroup with S-subgroups
$(L',\times_6 )$ and $(L'',\times_6 )$, $L'=\{2,4\}$ and
$L''=\{1,5\}$. This can be deduced from its multiplication table,
below. The triple $(U,V,W)$ such that
\begin{displaymath}
U=\left(\begin{array}{cccccc} 0 & 1 & 2 & 3 & 4 & 5\\
4 & 3 & 5 & 1 & 2 & 0
\end{array}\right),~V=
\left(\begin{array}{cccccc}
0 & 1 & 2 & 3 & 4 & 5 \\
1 & 3 & 2 & 4 & 5 & 0
\end{array}\right)~\textrm{and}~W=
\left(\begin{array}{cccccc}
0 & 1 & 2 & 3 & 4 & 5\\
1 & 0 & 5 & 4 & 2 & 3
\end{array}\right)
\end{displaymath}
are permutations on $L$, is an S-isotopism of
$(\mathbb{Z}_6,\times_6)$ unto an S-semigroup $(\mathbb{Z}_6,\ast)$
with S-subgroups $(L''',\ast )$ and $(L'''',\ast )$, $L'''=\{2,5\}$
and $L''''=\{0,3\}$ as shown in the second table below. Notice that
$A(L')=L'''$ and $A(L'')=L''''$ for all $A\in\{U,V,W\}$
and $U, V, W : L'\rightarrow L'''$ and $U, V, W : L''\rightarrow L''''$ are all bijcetions.\\
\begin{center}
\begin{tabular}{|c||c|c|c|c|c|c|}
\hline
$\times_6$ & 0 & 1 & 2 & 3 & 4 & 5\\
\hline \hline
0 & 0 & 0 & 0 & 0 & 0 & 0 \\
\hline
1 & 0 & 1 & 2 & 3 & 4 & 5 \\
\hline
2 & 0 & 2 & 4 & 0 & 2 & 4\\
\hline
3 & 0 & 3 & 0 & 3 & 0 & 3 \\
\hline
4 & 0 & 4 & 2 & 0 & 4 & 2 \\
\hline
5 & 0 & 5 & 4 & 3 & 2 & 1 \\
\hline
\end{tabular}
\qquad\qquad\qquad\qquad\qquad\qquad\qquad\qquad\begin{tabular}{|c||c|c|c|c|c|c|}
\hline
$\ast$ & 0 & 1 & 2 & 3 & 4 & 5\\
\hline \hline
0 & 0 & 1 & 2 & 3 & 4 & 5 \\
\hline
1 & 4 & 1 & 1 & 4 & 4 & 1 \\
\hline
2 & 5 & 1 & 5 & 2 & 1 & 2\\
\hline
3 & 3 & 1 & 5 & 0 & 4 & 2 \\
\hline
4 & 1 & 1 & 1 & 1 & 1 & 1 \\
\hline
5 & 2 & 1 & 2 & 5 & 1 & 5 \\
\hline
\end{tabular}
\end{center}
\end{myexam}
From Example~\ref{exam:1} and Example~\ref{exam:2}, it is very clear
that the study of of S-isotopy of two S-groupoids or S-quasigroups
or S-semigroups or S-loops is independent of the S-subsemigroup or
S-subgroup that is in consideration. All results in this paper are
true for any given S-subsemigroups or S-subgroups of two S-isotopic
S-groupoids or S-quasigroups or S-semigroups or S-loops. More
examples of S-isotopic S-groupoids can be constructed using
S-groupoids in \cite{van1}.

\begin{myrem}
Taking careful look at Definition~\ref{1:1} and comparing it with
[Definition~4.4.1,\cite{phd75}], it will be observed that the author
did not allow the component bijections $U$,$V$ and $W$ in $(U, V,
W)$ to act on the whole S-loop $L$ but only on the
S-subloop(S-subgroup) $L'$. We feel this is necessary to adjust here
so that the set $L-L'$ is not out of the study. Apart from this, our
adjustment here will allow the study of Smarandache isotopy to be
explorable. Therefore, the S-isotopism and S-isomorphism here are
clearly special types of relations(isotopism and isomorphism) on the
whole domain into the whole co-domain but those of Vasantha
Kandasamy \cite{phd75} only take care of the structure of the
elements in the S-subloop and not the S-loop. Nevertheless, we do
not fault her study for we think she defined them to apply them to
some life problems as an applied algebraist.
\end{myrem}

To every loop $(L,\cdot )$ with automorphism group $AUM(L,\cdot )$,
there corresponds another loop. Let the set $H=(L,\cdot )\times
AUM(L,\cdot )$. If we define '$\circ$' on $H$ such that $(\alpha,
x)\circ (\beta, y)=(\alpha\beta, x\beta\cdot y)$ for all $(\alpha,
x),(\beta, y)\in H$, then $H(L,\cdot )=(H,\circ)$ is a loop as shown
in Bruck \cite{phd82} and is called the Holomorph of $(L,\cdot )$.
Let $(L,\cdot )$ be an S-quasigroup(S-loop) with S-subgroup
$(L',\cdot )$. Define the Smarandache automorphism of $L$ to be the
set $SAUM(L)=SAUM(L,\cdot )=\{\alpha\in AUM(L)~:~\alpha~:~L'\to
L'\}$. It is easy to see that $SAUM(L)\le AUM(L)$. So, $SAUM(L)$
will be called the Smarandache automorphism group(SAG) of $L$. Now,
set $H_S=(L,\cdot )\times SAUM(L,\cdot )$. If we define '$\circ$' on
$H_S$ such that $(\alpha, x)\circ (\beta, y)=(\alpha\beta,
x\beta\cdot y)$ for all $(\alpha, x),(\beta, y)\in H_S$, then
$H_S(L,\cdot )=(H_S,\circ)$ is a S-quasigroup(S-loop) with
S-subgroup $(H',\circ )$ where $H'=L'\times SAUM(L)$ and thus will
be called the Smarandache Holomorph(SH) of $(L,\cdot )$.
\newpage
\section{Main Results}
\begin{myth}\label{1:4}
Let $U=(L,\oplus)$ and $V=(L,\otimes )$ be initial S-quasigroups
such that $SAUM(U)$ and $SAUM(V)$ are conjugates in $SSYM(L)$ i.e
there exists a $\psi\in SSYM(L)$ such that for any $\gamma\in
SAUM(V)$, $\gamma =\psi^{-1}\alpha\psi$ where $\alpha\in SAUM(U)$.
Then, $H_S(U)\succsim H_S(V)$ if and only if $x\delta\otimes y\gamma
=(x\beta\oplus y)\delta~\forall~x,y\in L,~\beta\in SAUM(U)$ and some
$\delta,\gamma\in SAUM(V)$. Hence:
\begin{enumerate}
\item $\gamma\in SAUM(U)$ if and only if $(I,\gamma ,\delta )\in
SAUT(V)$.
\item if $U$ is a initial S-loop, then;
\begin{enumerate}
\item ${\cal L}_{e\delta}\in SAUM(V)$.
\item $\beta\in SAUM(V)$ if and only if ${\cal R}_{e\gamma}\in
SAUM(V)$.
\end{enumerate}
where $e$ is the identity element in $U$ and ${\cal L}_x$, ${\cal
R}_x$ are respectively the left and right translations mappings of
$x\in V$.
\item if $\delta =I$, then $|SAUM(U)|=|SAUM(V)|=3$ and so
$SAUM(U)$ and $SAUM(V)$ are boolean groups.
\item if $\gamma =I$, then $|SAUM(U)|=|SAUM(V)|=1$.
\end{enumerate}
\end{myth}
{\bf Proof}\\
Let $H_S(L,\oplus )=(H_S,\circ )$ and $H_S(L,\otimes )=(H_S,\odot
)$. $H_S(U)\succsim H_S(V)$ if and only if there exists a bijection
$\phi ~:~H_S(U)\to H_S(V)$ such that $[(\alpha ,x)\circ (\beta
,y)]\phi =(\alpha,x)\phi\odot (\beta ,y)\phi$ and $(H',\oplus
)\overset{\phi}{\cong}(H'',\otimes )$ where $H'=L'\times SAUM(U)$
and $H''=L''\times SAUM(V)$, $(L',\oplus)$ and $(L'',\otimes )$ been
the initial S-subquasigroups of $U$ and $V$. Define $(\alpha,x)\phi
=(\psi^{-1}\alpha\psi ,x\psi^{-1}\alpha\psi )~\forall~(\alpha ,x)\in
(H_S,\circ )$ where $\psi\in SSYM(L)$.

$H_S(U)\cong H_S(V)\Leftrightarrow (\alpha\beta ,x\beta\oplus y)\phi
=(\psi^{-1}\alpha\psi ,x\psi^{-1}\alpha\psi )\odot
(\psi^{-1}\beta\psi ,y\psi^{-1}\beta\psi )\Leftrightarrow
(\psi^{-1}\alpha\beta\psi ,(x\beta\oplus y)\psi^{-1}\alpha\beta\psi
)=(\psi^{-1}\alpha\beta\psi ,x\psi^{-1}\alpha\beta\psi\otimes
y\psi^{-1}\beta\psi )\Leftrightarrow (x\beta\oplus
y)\psi^{-1}\alpha\beta\psi=x\psi^{-1}\alpha\beta\psi\otimes
y\psi^{-1}\beta\psi\Leftrightarrow x\delta\otimes y\gamma
=(x\beta\oplus y)\delta$ where $\delta =\psi^{-1}\alpha\beta\psi$,
$\gamma =\psi^{-1}\beta\psi$.

Note that, $\gamma{\cal L}_{x\delta}=L_{x\beta}\delta$ and
$\delta{\cal R}_{y\gamma}=\beta R_y\delta~\forall~x,y\in L$. So,
when $U$ is an S-loop, $\gamma{\cal L}_{e\delta}=\delta$ and
$\delta{\cal R}_{e\gamma}=\beta\delta$. These can easily be used to
prove the remaining part of the theorem.

\begin{myth}\label{1:6}
Let $\mathfrak{F}$ be any class of variety of S-quasigroups(loops).
Let $U=(L,\oplus)$ and $V=(L,\otimes )$ be initial
S-quasigroups(S-loops) that are S-isotopic under the triple of the
form $(\delta^{-1}\beta ,\gamma^{-1},\delta^{-1})$ for all $\beta\in
SAUM(U)$ and some $\delta,\gamma\in SAUM(V)$ such that their SAGs
are non-trivial and are conjugates in $SSYM(L)$ i.e there exists a
$\psi\in SSYM(L)$ such that for any $\gamma\in SAUM(V)$, $\gamma
=\psi^{-1}\alpha\psi$ where $\alpha\in SAUM(U)$. Then, $U\in
\mathfrak{F}$ if and only if $V\in \mathfrak{F}$.
\end{myth}
{\bf Proof}\\
By Theorem~\ref{1:4}, $H_S(U)\cong H_S(V)$. Let $U\in \mathfrak{F}$,
then since $H(U)$ has an initial S-subquasigroup(S-subloop) that is
isomorphic to $U$ and that initial S-subquasigroup(S-subloop) is
isomorphic to an S-subquasigroup(S-subloop) of $H(V)$ which is
isomorphic to $V$, $V\in \mathfrak{F}$. The proof for the converse
is similar.


\begin{thebibliography}{99}
\bibitem{phd79} J. O. Adeniran (2005), {\it On holomorphic theory of a
class of left Bol loops}, Al.I.Cuza 51, 1, 23-28
\bibitem{phd30} R. Artzy (1959), {\it Crossed inverse and related
loops}, Trans. Amer. Math. Soc. 91, 3, 480--492.
\bibitem{phd149} A. S. Basarab (1967), {\it A class of WIP-loops},
Mat. Issled. 2(2), 3-24.
\bibitem{phd146} A. S. Basarab(1970), {\it Isotopy of WIP
loops}, Mat. Issled. 5, 2(16), 3-12.
\bibitem{phd148} A. S. Basarab (1973), {\it The Osborn loop}, Studies in the theory of quasigroups and loops, 193. Shtiintsa,
Kishinev, 12--18.
\bibitem{phd147} A. S. Basarab and A.
I. Belioglo (1979), {\it Universal automorphic inverse G-loops},
Quasigroups and loops, Mat. Issled. 71, 3--7.
\bibitem{phd170} A. S. Basarab and A.
I. Belioglo (1979), {\it UAI Osborn loops}, Quasigroups and loops,
Mat. Issled. 51, 8--16.
\bibitem{phd41} R. H. Bruck (1966), {\it A survey of binary systems}, Springer-Verlag, Berlin-G\"ottingen-Heidelberg, 185pp.
\bibitem{phd82} R. H. Bruck (1944), {\it Contributions to the theory of loops}, Trans. Amer. Math. Soc. 55, 245--354.
\bibitem{phd40} R. H. Bruck and L. J. Paige (1956), {\it Loops whose
inner mappings are automorphisms}, The annuals of Mathematics, 63,
2, 308--323.
\bibitem{phd39} O. Chein, H. O. Pflugfelder and J. D. H. Smith (1990), {\it Quasigroups and loops : Theory and applications}, Heldermann Verlag, 568pp.
\bibitem{phd80} V. O. Chiboka and A. R. T. Solarin (1991), {\it Holomorphs of conjugacy closed loops}, Scientific Annals of Al.I.Cuza. Univ. 37, 3, 277--284.
\bibitem{phd49} J. D\'{e}ne and A. D. Keedwell (1974), {\it Latin squares and their applications}, the English University press Lts, 549pp.
\bibitem{phd159} E. Falconer (1969), {\it Quasigroup identities invariant under
isotopy}, Ph.D thesis, Emory University.
\bibitem{phd160} E. Falconer (1970), {\it Isotopy invariants in
quasigroups}, Trans. Amer. Math. Soc. 151, 2, 511-526.
\bibitem{phd50} F. Fenyves (1968), {\it Extra loops I}, Publ. Math. Debrecen, 15, 235--238.
\bibitem{phd56} F. Fenyves (1969), {\it Extra loops II}, Publ. Math. Debrecen, 16, 187--192.
\bibitem{phd42} E. G. Goodaire, E. Jespers and C. P. Milies (1996), {\it Alternative loop rings}, NHMS(184), Elsevier, 387pp.
\bibitem{phd44} E. D. Huthnance Jr.(1968), {\it A theory of
generalised Moufang loops}, Ph.D. thesis, Georgia Institute of
Technology.
\bibitem{phd143} T. G. Ja\'iy\'e\d ol\'a (2005), {\it An isotopic study of
properties of central loops}, M.Sc. dissertation, University of
Agriculture, Abeokuta.
\bibitem{sma1} T. G. Ja\'iy\'e\d ol\'a (2006), {\it An holomorphic study of the Smarandache concept in
loops}, Scientia Magna Journal, 2, 1, 1--8.
\bibitem{sma2} T. G. Ja\'iy\'e\d ol\'a (2006), {\it Parastrophic invariance of Smarandache quasigroups},
Scientia Magna Journal, 2, 3, 48--53.
\bibitem{sma3} T. G. Ja\'iy\'e\d ol\'a (2006), {\it On the universality of some Smarandache loops of Bol-
Moufang type}, Scientia Magna Journal, 2, 4, 45-48.
\bibitem{phd150} B. B. Karklinsh and V. B. Klin (1976), {\it
Universal automorphic inverse property loops}, Collections: sets and
quasigroups, Shtintsa.
\bibitem{phd95} T.
Kepka, M. K. Kinyon, J. D. Phillips, {\it The structure of
F-quasigroups}, communicated for publication.
\bibitem{phd118} T. Kepka, M. K. Kinyon, J. D. Phillips, {\it
F-quasigroups and generalised modules}, communicated for
publication.
\bibitem{phd119} T. Kepka, M. K. Kinyon, J. D. Phillips, {\it F-quasigroups isotopic to groups}, communicated for publication.
\bibitem{phd33} M. K. Kinyon (2005), {\it A survey of Osborn loops},
Milehigh conference on loops, quasigroups and non-associative
systems, University of Denver, Denver, Colorado.
\bibitem{phd124} M. K. Kinyon, J. D. Phillips and P. Vojt\v echovsk\'y (2005), {\it Loops of Bol-Moufang type with a subgroup of index
two}, Bul. Acad. Stiinte Repub. Mold. Mat. 3(49), 71--87.
\bibitem{muk1} A. S. Muktibodh (2005), {\it Smarandache quasigroups rings},
Scientia Magna Journal, 1, 2, 139--144.
\bibitem{muk2} A. S. Muktibodh (2006), {\it Smarandache Quasigroups},
Scientia Magna Journal, 2, 1, 13--19.
\bibitem{phd89} J. M. Osborn (1961), {\it Loops with the weak
inverse property}, Pac. J. Math. 10, 295--304.
\bibitem{phd141} Y. T. Oyebo and O. J. Adeniran, {\it On the holomorph of central
loops}, Pre-print.
\bibitem{phd3} H. O. Pflugfelder (1990), {\it Quasigroups and loops : Introduction}, Sigma series in Pure Math. 7, Heldermann Verlag, Berlin, 147pp.
\bibitem{phd9} J. D. Phillips and P. Vojt\v echovsk\'y (2005), {\it The varieties of loops of Bol-Moufang type}, Alg. Univer. 3(54), 259--383..
\bibitem{phd61} J. D. Phillips and P. Vojt\v echovsk\'y (2005), {\it The varieties of quasigroups of Bol-Moufang type : An equational
approach}, J. Alg. 293, 17--33.
\bibitem{phd85} D. A. Robinson (1964), {\it Bol loops}, Ph. D thesis,
University of Wisconsin, Madison, Wisconsin.
\bibitem{phd7} D. A. Robinson (1971), {\it Holomorphic theory of
extra loops}, Publ. Math. Debrecen 18, 59--64.
\bibitem{phd97} P. N. Syrbu (1996), {\it On loops with universal elasticity}, Quasigroups and Related
Systems, 3, 41--54.
\bibitem{phd75} W. B. Vasantha Kandasamy (2002), {\it Smarandache
loops}, Department of Mathematics, Indian Institute of Technology,
Madras, India, 128pp.
\bibitem{phd83} W. B. Vasantha Kandasamy (2002), {\it Smarandache
Loops}, Smarandache Notions Journal, 13, 252--258.
\bibitem{van1} W. B. Vasantha Kandasamy (2002), {\it Groupoids and Smarandache Groupoids},
American Research Press Rehoboth, 114pp.
\bibitem{van2} W. B. Vasantha Kandasamy (2002), {\it Smarandache Semigroups},
American Research Press Rehoboth, 94pp.
\bibitem{van3} W. B. Vasantha Kandasamy (2002), {\it Smarandache Semirings, Semifields, And Semivector Spaces},
American Research Press Rehoboth, 121pp.
\bibitem{van4} W. B. Vasantha Kandasamy (2003), {\it Linear Algebra And Smarandache Linear Algebra},
American Research Press, 174pp.
\bibitem{van5} W. B. Vasantha Kandasamy (2003), {\it Bialgebraic Structures And Smarandache Bialgebraic Structures},
American Research Press Rehoboth, 271pp.
\bibitem{van6} W. B. Vasantha Kandasamy and F. Smarandache (2005), {\it N-Algebraic Structures And Smarandache N-Algebraic Structures},
Hexis Phoenix, Arizona, 174pp.
\end{thebibliography}
\end{document}